# THE NON-AMENABILITY OF SCHREIER GRAPHS FOR INFINITE INDEX QUASICONVEX SUBGROUPS OF HYPERBOLIC GROUPS

ILYA KAPOVICH

ABSTRACT. We show that if $H$ is a quasiconvex subgroup of infinite index in a non-elementary hyperbolic group $G$ then the Schreier coset graph for $G$ relative to $H$ is non-amenable.

## 1. INTRODUCTION

A connected graph of bounded degree $X$ is *non-amenable* if $X$ has nonzero Cheeger constant, or, equivalently, if the spectral radius of the simple random walk on $X$ is less than one (see Section 2 below for more precise definitions). Non-amenable graphs play an increasingly important role in the study of various probabilistic phenomena, such as random walks, harmonic analysis, Brownian motion and percolations, on graphs and manifolds (see for example [2, 5, 6, 7, 14, 16, 17, 23, 29, 42, 43, 60, 69, 70]) as well as in the study of expander families of finite graphs (see for example [50, 64, 65]).

It is well-known that a finitely generated group $G$ is non-amenable if and only if some (any) Cayley graph of $G$ with respect to a finite generating set is non-amenable. Word-hyperbolic groups are non-amenable unless they are virtually cyclic and thus their Cayley graphs provide a large and interesting class of non-amenable graphs. In this paper we investigate non-amenability of Schreier coset graphs corresponding to subgroups of hyperbolic groups.

We recall the definition of a Schreier coset graph:

**Definition 1.1.** Let $G$ be a group and let $\pi : A \to G$ be a map where $A$ is a finite alphabet such that $\pi(A)$ generates $G$ (we refer to such $A$ as a *marked finite generating set* or just *finite generating set* of $G$). Let $H \le G$ be a subgroup of $G$. The *Schreier coset graph* (or the *relative Cayley graph*) $\Gamma(G, H, A)$ for $G$ relative to $H$ with respect to $A$ is an oriented labeled graph defined as follows:

1. The vertices of $\Gamma = \Gamma(G, H, A)$ are precisely the cosets of $H$ in $G$, that is $V\Gamma := \{Hg \,|\, g \in G\}$.
2. The set of positively oriented edges of $\Gamma(G, H, A)$ is in one-to-one correspondence with the set $V\Gamma \times A$. For each pair $(Hg, a) \in V\Gamma \times A$ there is a positively oriented edge in $\Gamma(G, H, A)$ from $Hg$ to $Hg\pi(a)$ labeled by the letter $a$.

Thus the label of every path in $\Gamma(G, H, A)$ is a word in the alphabet $A \cup A^{-1}$. The graph $\Gamma(G, H, A)$ is connected since $\pi(A)$ generates $G$. Moreover, $\Gamma(G, H, A)$ comes equipped with the natural simplicial metric $d$ obtained by giving every edge length one.

We can identify the Schreier graph with 1-skeleton of the presentation complex of $G$ corresponding to any presentation of $G$ of the form $G = \langle A \,|\, R \rangle$. It is also easy to see that if $M$ is a closed Riemannian manifold and $H \leq G = \pi_1(M)$, then the Schreier graph $\Gamma(G, H, A)$ is quasi-isometric to the covering space of $M$ corresponding to $H$.

If $H$ is normal in $G$ and $G_1 = G/H$ is the quotient group, then $\Gamma(G, H, A)$ is exactly the Cayley graph of the group $G_1$ with respect to $A$. In particular, if $H = 1$ then $\Gamma(G, 1, A)$ is the standard *Cayley graph of $G$ with respect to $A$*, denoted $\Gamma(G, A)$.

The notion of a word-hyperbolic group was introduced by M.Gromov [39]. Recall that a finitely generated group $G$ is said to be *word-hyperbolic* if for any finite generating set $A$ of $G$ there is $\delta \geq 0$ such that all geodesic triangles in $\Gamma(G, A)$ are $\delta$-thin, that is each side of a triangle is contained in the closed $\delta$-neighborhood of the union of the other two sides. A subgroup $H$ of a word-hyperbolic group $G$ is said to be *quasiconvex* in $G$ if for any finite generating set $A$ of $G$ there is $\epsilon \geq 0$ such that every geodesic in $\Gamma(G, A)$ with both endpoints in $H$ is contained in the $\epsilon$-neighborhood of $H$ in $\Gamma(G, A)$. Quasiconvex subgroups are closely related to geometric finiteness in the Kleinian group context [67]. They enjoy a number of particularly good properties and play an important role in hyperbolic group theory and its applications (see for example [3, 4, 8, 30, 33, 34, 35, 36, 37, 41, 44, 45, 47, 49, 51, 53, 59, 68]).

Our main result is the following:

**Theorem 1.2.** *Let $G$ be a non-elementary word-hyperbolic group with a marked finite generating set $A$. Let $H \leq G$ be a quasiconvex subgroup of infinite index in $G$. Then the Schreier coset graph $\Gamma(G, H, A)$ is non-amenable.*

The study of Schreier graphs arises naturally in various generalizations of J.Stallings' theory of ends of groups [22, 28, 58, 59, 61]. The case of virtually cyclic (and hence quasiconvex) subgroups of hyperbolic groups is particularly important to understand in the theory of JSJ-decomposition for hyperbolic groups originally developed by Z.Sela [63] and later by B.Bowditch [9] (see also [57, 22, 27, 62] for various generalizations of the JSJ-theory). A variation of the Følner criterion of non-amenability (see Proposition 2.3 below), when the Cheeger constant is defined by taking the infimum over all finite subsets containing no more than a half of all the vertices, is used to define an important notion of *expander families* of finite graphs. Most known sources of expander families involve taking Schreier coset graphs corresponding to subgroups of finite index in a group with Kazhdan property (T) (see [50, 64, 65] for a detailed exposition on expander families and their connections with non-amenability).

Since non-amenable graphs of bounded degree are well-known to be transient with respect to the simple random walk, Theorem 1.2 implies that $\Gamma(G, H, A)$ is also transient. This was first shown in [48] by more elementary means for the case when $G$ is torsion-free hyperbolic $H \leq G$ is quasiconvex of infinite index.

It was originally stated by M.Gromov [39] and proved by R.Foord [26] and I.Kapovich [48] that for any quasiconvex subgroup $H$ in a hyperbolic group $G$ with a finite generating set $A$ the coset graph $\Gamma(G, H, A)$ is a hyperbolic metric space. A great deal is known about random walks on hyperbolic graphs, but most of these results assume some kind of non-amenability. Thus Theorem 1.2 together with hyperbolicity of $\Gamma(G, H, A)$ and a result of A.Ancona [2] (see also the [70]) immediately imply:

**Corollary 1.3.** *Let $G$ be a non-elementary word-hyperbolic group with a finite generating set $A$. Let $H \leq G$ be a quasiconvex subgroup of infinite index in $G$ and let $Y$ be the Schreier coset graph $\Gamma(G, H, A)$.*
*Then:*

1. *The trajectory of almost every simple random walk on $Y$ converges in the topology of $Y \cup \partial Y$ to some point in $\partial Y$ (where $\partial Y$ is the hyperbolic boundary).*
2. *There Martin boundary of a the simple random walk on $X$ is homeomorphic to the hyperbolic boundary $\partial X$ and the Martin compactification $\hat{X}$ for the simple random walk on $X$ is homeomorphic to the hyperbolic compactification $X \cup \partial X$.*

The statement of Theorem 1.2 is easy to illustrate for the case of a free group. Suppose $F = F(a, b)$ is free and $H \leq F$ is a finitely generated subgroup of infinite index (which is therefore quasiconvex [66]). Put $A = \{a, b\}$. Then the Schreier graph $Y = \Gamma(F, H, A)$ looks like a finite graph with several infinite tree-brunches attached to it (the "brunches" are 4-regular trees except for the attaching vertices). In this situation it is easy to see that $Y$ has positive Cheeger constant and so $Y$ is non-amenable. Alex Lubotzky and Andrzej Zuk pointed out to the author that if $G$ is a group with Kazhdan property (T) then for any subgroup $H$ of infinite index in $G$ the Schreier coset graph for $G$ relative to $H$ is non-amenable. There are many examples of word-hyperbolic groups with Kazhdan property (T) and in view of Theorem 1.2 it would be particularly interesting to investigate if they can possess non-quasiconvex finitely generated subgroups.

Non-amenability of graphs is closely related to co-growth. Thus we also obtain the following fact.

**Corollary 1.4.** *Let $G = \langle x_1, \ldots, x_k \, | \, r_1, \ldots, r_m \rangle$ be a non-elementary hyperbolic group and let $H \leq G$ be a quasiconvex subgroup of infinite index. Let $a_n$ be the number of freely reduced words in $A = \{x_1, \ldots, x_k\}^{\pm 1}$ of length $n$ representing elements of $H$. Let $b_n$ be the number of all words in $A$ of length $n$ representing elements of $H$.*

*Then*
$$\limsup_{n\to\infty} \sqrt[n]{a_n} < 2k-1$$

*and*
$$\limsup_{n\to\infty} \sqrt[n]{b_n} < 2k.$$

It is easy to see that the statement of Theorem 1.2 need not hold for finitely generated subgroups which are not quasiconvex. For example, a celebrated construction of E.Rips [56] states that for any finitely presented group $Q$ there is a short exact sequence

$$1 \to K \to G \to Q \to 1,$$

where $G$ is non-elementary torsion-free word-hyperbolic and where $K$ can be generated by two elements (but $K$ is not finitely presentable). If $Q$ is chosen to be non-amenable, then the Schreier graph for $G$ relative to $H$ is non-amenable. Finitely presentable and even hyperbolic examples are also possible. For instance, if $F$ is a free group of finite rank and $\phi : F \to F$ is an atoroidal automorphism, then the mapping torus group of $\phi$

$$M_\phi = \langle F, t \,|\, t^{-1}ft = \phi(f) \text{ for all } f \in F \rangle$$

is word-hyperbolic [8, 11]. In this case $G/F \simeq \mathbb{Z}$ and thus amenable.

As a main technical tool in the proof of Theorem 1.2 we obtain the following result, which appears to be of independent interest.

**Theorem 1.5.** *Let $G$ be a non-elementary hyperbolic group and let $H \leq G$ be a quasiconvex subgroup of infinite index. Then there exists a free quasiconvex subgroup of rank two $F \leq G$ which is conjugacy separated from $H$ in $G$, that is for any $g \in G$*

$$g^{-1}Hg \cap F = 1.$$

The author is grateful to Laurent Bartholdi, Philip Bowers, Christophe Pittet and Tatiana Smirnova-Nagnibeda for the many helpful discussions regarding random walks and to Paul Schupp for encouragement.

## 2. Non-amenability for graphs

Let $X$ be a connected graph of bounded degree. We will denote by $\rho(X)$ the *spectral radius* of $X$ which can be defined as

$$\rho(X) := \limsup_{n\to\infty} \sqrt[n]{p^{(n)}(x,y)}$$

where $x, y$ are two vertices of $X$ and $p^{(n)}(x,y)$ is the probability that a simple random walk starting at $x$ will end up at $y$ in $n$ steps. It is well-known that $\rho(X) \leq 1$ and that the definition of $\rho(X)$ does not depend on the choice of $x, y$.

**Definition 2.1** (Amenability for graphs). A connected graph $X$ of bounded degree is said to be *amenable* if $\rho(X) = 1$ and *non-amenable* if $\rho(X) < 1$.

It is also well-known that non-amenability of $X$ implies that $X$ is transient (see for example Theorem 51 of [15]). We refer the reader to [15, 69, 70] for the comprehensive background information about random walks on graphs and for further references on this topic.

**Convention 2.2.** Let $X$ be a connected graph of bounded degree with the simplicial metric $d$. For a finite nonempty subset $S \subset VX$ we will denote by $|S|$ the number of elements in $S$.

If $S$ is a finite subset of the vertex set of $X$ and $k \geq 1$ is an integer, we will denote by $\mathcal{N}_k^X(S) = \mathcal{N}_k(S)$ the set of all vertices $v$ of $X$ such that $d_X(v, S) \leq k$. Also, we will denote $\eth^X S = \eth S := \mathcal{N}_1(S) - S$.

The number

$$\iota(X) := \inf\{\frac{|\eth S|}{|S|} \: : \: S \text{ is a finite nonempty subset of the vertex set of } X\}$$

is called the *Cheeger constant* or the *isoperimetric constant* of $X$.

There are many alternative definitions of non-amenability:

**Proposition 2.3.** *Let $X$ be a connected graph of bounded degree with simplicial metric $d$. Then the following conditions are equivalent:*

1. *The graph $X$ is non-amenable.*
2. *(Følner criterion) We have $\iota(X) > 0$.*
3. *(Gromov's Doubling Condition) There is some $k \geq 1$ such that for any finite nonempty subset $S \subseteq VX$ we have*

   $$|\mathcal{N}_k(S)| \geq 2|S|.$$

4. *For any integer $q > 1$ there is some $k \geq 1$ such that for any finite nonempty subset $S \subseteq VX$ we have*

   $$|\mathcal{N}_k(S)| \geq q|S|.$$

5. *For some $0 < \sigma < 0$ $p^{(n)}(x, y) = o(\sigma^n)$ for any $x, y \in VX$.*
6. *The pseudogroup $W(X)$ consisting of all bijections between subsets of $VX$ which are "bounded perturbations of the identity" admits a "paradoxical decomposition" (see [15] for definitions).*
7. *("Grasshoper criterion") There exists a map $\phi : VX \to VX$ such that $\sup_{x \in VX} d(x, \phi(x)) < \infty$ and that for any $x \in VX$ $|\phi^{-1}(x)| \geq 2$.*
8. *There exists a map $\phi : VX \to VX$ such that $\sup_{x \in VX} d(x, \phi(x)) < \infty$ and that for any $x \in VX$ $|\phi^{-1}(x)| = 2$.*
9. *The bottom of the spectrum for the combinatorial Laplacian operator on $X$ is $> 0$ (see [20] for the precise definitions).*
10. *We have $H_0^{uf}(X) = 0$ (see [12] for the precise definition of* uniformly finite *homology groups $H_i^{uf}$).*

11. We have $H_0^{(l_p)}(X) = 0$ for any $1 < p < \infty$ (see [23] for the precise definition of $H_i^{(l_p)}$).

*Proof.* All of the above statements are well-known, but we will still provide some sample references.

The fact that (1), (2), (5) and (6) are equivalent is stated in Theorem 51 of [15]. The fact that (3), (4), (6), (7) and (8) are equivalent follows from Theorem 32 of [15]. The equivalence of (2) and (9) is due to J.Dodziuk [20]. J.Block and S.Weinberger [12] established the equivalence of (2) and (10). Finally, G.Elek [23] proved that (2) is equivalent to (11). □

In case of regular graphs one can also characterize non-amenability in terms of co-growth.

**Definition 2.4.** Let $X$ be a connected graph of bounded degree with a base-vertex $x_0$. Let $a_n = a_n(X, x_0)$ be the number of reduced edge-paths of length $n$ from $x_0$ to $x_0$. Let $b_n = b_n(X, x_0)$ be the number of all edge-paths of length $n$ from $x_0$ to $x_0$. Put

$$\alpha(X) := \limsup_{n \to \infty} \sqrt[n]{a_n} \text{ and } \alpha(X) := \limsup_{n \to \infty} \sqrt[n]{b_n}$$

Then we will call $\alpha(X)$ the *co-growth rate* of $X$ and we will call $\beta(X)$ the *non-reduced co-growth rate* of $X$. These definitions are well-known to be independent of the choice of $x_0$.

It is easy to see that for a $d$-regular connected graph $X$ we have $\alpha(X) \leq d - 1$ and $\beta(X) \leq d$. Moreover, $\rho(X) = \frac{\beta(X)}{d}$. It turns out that non-amenability of regular graphs can be characterized in terms of the co-growth rate. The following result is was originally proved by R.Grigorchuk [38] and J.Cohen [18] for Cayley graphs of finitely generated groups and by L.Bartholdi [5] for arbitrary regular graphs.

**Theorem 2.5.** [5] *Let $X$ be a connected $d$-regular graph with $d \geq 3$. Put $\alpha = \alpha(X)$, $\beta = \beta(X)$ and $\rho = \rho(X)$. Then*

$$\rho = \frac{2\sqrt{d-1}}{d} \quad \text{if} \quad 1 \leq \alpha \leq \sqrt{d-1}$$

*and*

$$\rho = \frac{\sqrt{d-1}}{d}\left(\frac{\sqrt{d-1}}{\alpha} + \frac{\alpha}{\sqrt{d-1}}\right) \quad \text{if} \quad \sqrt{d-1} \leq \alpha \leq d-1.$$

*In particular* $\rho < 1 \iff \alpha < d-1 \iff \beta < d$.

## 3. Hyperbolic metric spaces

The basic information about Gromov-hyperbolic metric spaces and word-hyperbolic groups can be found in [39, 19, 31, 1, 13, 24, 4] and other sources. We will briefly recall the main definitions.

If $(X, d)$ is a geodesic metric space and $x, y \in X$, we shall denote by $[x, y]$ a geodesic segment from $x$ to $y$ in $X$. Also, of $p$ is a path in $X$, we will denote the length of $p$ by $l(p)$. Two paths $\alpha$ and $\beta$ in $X$ are said to be *K-Hausdorff close* if each of them is contained in the closed $K$-neighborhood of the other. Two paths $\alpha, \beta : [0, T] \to X$ (where $T \geq 0$) are said to be *K-uniformly close* (or to be *K-fellow travelers*) if for any $0 \leq t \leq T$ we have $d(\alpha(t), \beta(t)) \leq K$.

Given a path $\alpha : [0, T] \to X$ we shall often identify $\alpha$ with its image $\alpha([0, T]) \subseteq X$. If $Z \subseteq X$ and $\epsilon \geq 0$, we will denote the closed $\epsilon$-neighborhood of $Z$ in $X$ by $N_\epsilon(Z)$.

**Definition 3.1** (Hyperbolic metric space). Let $(X, d)$ be a geodesic metric space and let $\delta \geq 0$. The space $X$ is said to be $\delta$-*hyperbolic* if for any geodesic triangle in $X$ with sides $\alpha, \beta, \gamma$ we have

$$\alpha \subseteq N_\delta(\beta \cup \gamma), \ \beta \subseteq N_\delta(\alpha \cup \gamma) \text{ and } \gamma \subseteq N_\delta(\alpha \cup \beta)$$

that is for any $p \in \alpha$ there is $q \in \beta \cup \gamma$ such that $d(p, q) \leq \delta$ (and the symmetric condition holds for any $p \in \beta$ and any $p \in \gamma$).

A geodesic space $X$ is said to be *hyperbolic* if it is $\delta$-hyperbolic for some $\delta \geq 0$.

Suppose $\alpha$ and $\beta$ are geodesic segments from $x$ to $y$ in a geodesic metric space $(X, d)$. We will say that the geodesic bigon $\Theta = \alpha \cup \beta$ is $\delta$-*thin* if $\alpha \subseteq N_\delta(\beta)$ and $\beta \subseteq N_\delta(\alpha)$.

We will also need another equivalent definition of hyperbolicity.

**Definition 3.2** (Gromov product). Let $(X, d)$ be a metric space and suppose $x, y, z \in X$. We set

$$(x, y)_z := \frac{1}{2}[d(z, x) + d(z, y) - d(x, y)]$$

Note that $(x, y)_z = (y, x)_z$.

**Definition 3.3.** Let $(X, d)$ be a geodesic metric space and let $\Delta = \alpha_1 \cup \alpha_2 \cup \alpha_3$ is a geodesic triangle in $X$, where $\alpha_1 = [z, x]$, $\alpha_2 = [z, y]$ and $\alpha_3 = [x, y]$.

Note that by definition of Gromov product $d(x, y) = (z, y)_x + (x, z)_y$, $d(x, z) = (y, z)_x + (x, y)_z$ and $d(y, z) = (x, z)_y + (x, y)_z$. Thus there exist uniquely defined points $p \in \alpha_1$, $q \in \alpha_2$ and $r \in \alpha_3$ such that:

$$d(z, p) = d(z, q) = (x, y)_z, d(x, p) = d(x, r) = (y, z)_x, \quad \text{and}$$
$$d(y, q) = d(y, r) = (x, z)_y$$

We will call $(p, q, r)$ the *inscribed triple* of $\Delta$.

**Definition 3.4** (Trim triangle). Let $(X, d)$ be a geodesic metric space and let $\delta \geq 0$. Let $\Delta = \alpha_1 \cup \alpha_2 \cup \alpha_3$ be a geodesic triangle in $X$, where $\alpha_1 = [z, x]$, $\alpha_2 = [z, y]$ and $\alpha_3 = [x, y]$. We say that $\Delta$ is $\delta$-*trim* if the following holds.

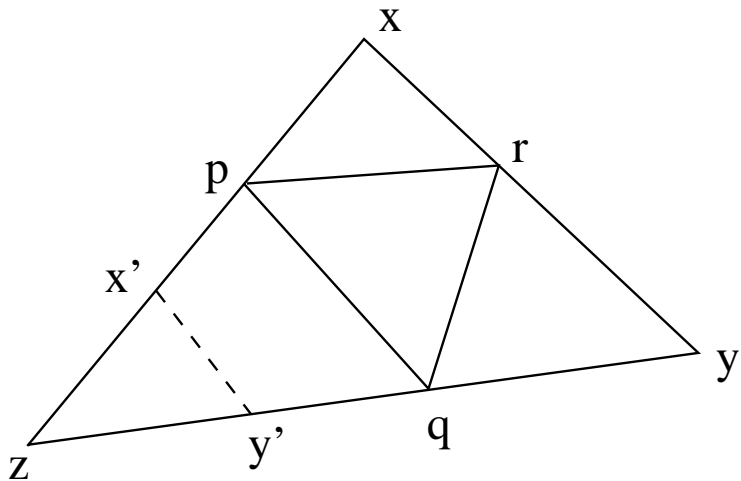

Figure 1. Trim triangle

Let $(p, q, r)$ be the inscribed triple of $\Delta$, where $p \in \alpha_1, q \in \alpha_2, r \in \alpha_3$, as shown in Figure 1.

Then:

1. for any points $x' \in \alpha_1, y' \in \alpha_2$ with $d(z, x') = d(z, y') \leq (x, y)_z$ we have $d(x', y') \leq \delta$;
2. for any points $z' \in \alpha_1, y' \in \alpha_3$ with $d(x, z') = d(x, y') \leq (z, y)_x$ we have $d(z', y') \leq \delta$;
3. for any points $z' \in \alpha_2, x' \in \alpha_3$ with $d(y, x') = d(y, z') \leq (x, z)_y$ we have $d(x', z') \leq \delta$.

The following statement is well-known [1]:

**Theorem 3.5.** *Let $(X, d)$ be a geodesic metric space. Then $(X, d)$ is hyperbolic if and only if for some $\delta \geq 0$ all geodesic triangles in $X$ are $\delta$-trim.*

We also shall make use of the following property of the Gromov product which in fact can be used as an equivalent definition of hyperbolicity.

**Proposition 3.6.** [1] *Let $(X, d)$ be a geodesic metric space with $\delta$-trim geodesic triangles, where $\delta \geq 0$. Then for any $p, x, y, z \in X$ we have:*

$$(x, y)_p \geq \min\{(x, z)_p, (y, z)_p\} - 2\delta.$$

It is well-known and easy to check from the definition that Gromov product is monotone in both arguments:

**Lemma 3.7.** *Let $(X, d)$ be a geodesic metric space. Suppose $[x, y]$ and $[x, z]$ are geodesic segments in $X$. Then for any $y' \in [x, y]$ and any $z' \in [x, z]$ we have $(y', z')_x \leq (y, z)_x$.*

**Definition 3.8** (Quasigeodesics and local quasigeodesics)**.** Suppose $(X, d)$ is a geodesic metric space. For $\lambda \geq 1, \epsilon \geq 0$ we will say that a path $\alpha : [a, b] \to X$ is a $(\lambda, \epsilon)$-*quasigeodesic*, if $\alpha$ is parameterized by arc-length and for any $s_1, s_2 \in [a, b]$ we have $|s_2 - s_1| \leq \lambda d(\alpha(s_1), \alpha(s_1)) + \epsilon$. A path $\alpha$ in $X$ is said to be $\lambda$-*quasigeodesic*, if it is $(\lambda, \lambda)$-quasigeodesic. A naturally parameterized path $\alpha : [a, b] \to X$ is a $N$-*local $\lambda$-quasigeodesic* if for any $a \leq s_1 \leq s_2 \leq b$ with $|s_2 - s_1| \leq N$ the restriction $\alpha|_{[s_1, s_2]}$ is a $\lambda$-quasigeodesic.

It is well-known that in hyperbolic spaces local quasigeodesics are global quasigeodesics, provided the local parameter is sufficiently big:

**Proposition 3.9** (Pasting quasigeodesics)**.** [19, 31] *For any $\delta \geq 0$ and $\lambda \geq 1$ there exist constants $N = N(\delta, \lambda) > 0$ and $\lambda' = \lambda'(\delta, \lambda) \geq 1$ with the following property.*

*Suppose $(X, d)$ is a $\delta$-hyperbolic geodesic metric space. Then any $N$-local $\lambda$-quasigeodesic in $X$ is a $\lambda'$-quasigeodesic.*

**Definition 3.10** (Gromov product for sets)**.** Let $(X, d)$ be a metric space. Let $x \in X$ and $Q, Q' \subseteq X$. Put $(Q, Q')_x := sup\{(q, q')_x \,|\, q \in Q, q' \in Q'\}$.

**Lemma 3.11.** *Let $(X, d)$ be a geodesic metric space with $\delta$-trim triangles. Suppose $x \in X$ and $Q, Q', Q''$ are nonempty subsets of $X$ such that $(Q', Q)_x > L$, $(Q'', Q)_x > L$. Then $(Q', Q'')_x \geq L - 2\delta$.*

*Proof.* This statement follows immediately from Proposition 3.6. □

## 4. Quasiconvex subgroups of hyperbolic groups

The detailed background information on quasiconvex subgroups of hyperbolic groups can be found in [1, 19, 31, 66, 4, 30, 49, 52, 37, 33] and other sources.

**Convention 4.1.** Let $G$ be a fixed non-elementary hyperbolic group with a fixed finite generating set $A$. Let $X = \Gamma(G, A)$ be the Cayley graph of $G$ with respect to $A$. We will denote the word-metric corresponding to $A$ on $X$ by $d$. Also, for $g \in G$ we will denote $|g|_A := d_A(1, g)$. Let $\delta \geq 10$ be an integer such that $(X, d)$ has $\delta$-trim geodesic triangles. For an $A$-word $w$ we will denote by $l(w)$ the length of $w$ and by $\overline{w}$ the element of $G$ represented by $w$. A word $w$ in $A \cup A^{-1}$ will be called an *geodesic* if $|w| = |\overline{w}|_A$, that is if any path labeled $w$ in $X$ is geodesic. Similarly, a word $w$ is said to be $\lambda$-*quasigeodesic* if any path labeled $w$ in $X$ is $\lambda$-quasigeodesic. If $s \subseteq G$, we will denote by $\langle S \rangle$ the subgroup of $G$ generated by $G$. For a word $w$ we will denote the length of $w$ by $|w|$.

These constants, notations and conventions will be fixed till the end of this article, unless specified otherwise.

**Definition 4.2** (Quasiconvexity)**.** A subset $Z \subseteq X$ is called $\epsilon$-*quasiconvex*, where $\epsilon \geq 0$, if for any $z_1, z_2 \in Z$ and any geodesic $[z_1, z_2]$ in $X$ we have $N_\epsilon(Z)$. A subset $Z \subseteq X$ is *quasiconvex* if it is $\epsilon$-quasiconvex for some $\epsilon \geq 0$. A subgroup $H \leq G$ is *quasiconvex* if $H \subseteq X$ is a quasiconvex subset.

We summarize some well-known basic facts regarding quasiconvex subgroups:

**Proposition 4.3.** *Let $G$ be a word-hyperbolic group with a finite generating set $A$. Let $X = \Gamma(G, A)$ be the Cayley graph of $G$ with the word-metric $d$ induced by $A$. Let $\delta \geq 1$ be an integer such that all geodesic triangles in $X$ are $\delta$-trim. Then:*

1. *[19, 31] If $H \leq G$ is a subgroup, then either $H$ is virtually cyclic (in which case $H$ is called* elementary*) or $H$ contains a free subgroup $F$ of rank two which is quasiconvex in $G$ (in which case $H$ is said to be* non-elementary*).*
2. *[19, 31, 4, 30] Let $B$ be another finite generating set of $G$ and let $Y = \Gamma(G, B)$. Suppose $H \leq G$ is a subgroup. Then $H \subseteq X$ is quasiconvex if and only if $H \subseteq Y$ is quasiconvex.*
3. *[1, 19, 31] Every cyclic subgroup of $G$ is quasiconvex in $G$.*
4. *[1, 19, 31] If $H \leq G$ is quasiconvex then $H$ is finitely presentable and word-hyperbolic.*
5. *[19, 31, 4, 30] Suppose $H \leq G$ is generated by a finite set $Q$ inducing a word-metric $d_Q$ on $H$. Then $H$ is quasiconvex in $G$ if and only if there is $C > 0$ such that for any $h_1, h_2 \in H$*
$$d_Q(h_1, h_2) \leq Cd(h_1, h_2).$$
6. *[30] The set $\mathcal{L}$ of all $A$-geodesic words is a regular language which provides a bi-automatic structure for $G$. Moreover, a subgroup $H \leq G$ is quasiconvex if and only if $H$ is $\mathcal{L}$-rational, that is the set $\mathcal{L}_H = \{w \in \mathcal{L} \mid \overline{w} \in H\}$ is a regular language.*
7. *[66] If $H_1, H_2 \leq G$ are quasiconvex, then $H_1 \cap H_2 \leq G$ is quasiconvex.*
8. *[49] Suppose $H \leq G$ is an infinite quasiconvex subgroup. Then $H$ has finite index in its* commensurator $Comm_G(H)$, *where*

$Comm_G(H) :=$
$= \{g \in G \mid [H : g^{-1}Hg \cap H] < \infty \text{ and } [g^{-1}Hg : Hg^{-1}Hg \cap H] < \infty\}.$

The following useful fact follows directly from the proofs of Lemma 4.1 and Lemma 4.5 of [4]:

**Lemma 4.4.** *For any quasiconvex subgroup $H \leq G$ there exists an integer constant $K = K(G, H, A) > 0$ with the following properties.*

*Suppose $g \in G$ is shortest with respect to $d$ in the coset class $Hg$. Let $h \in H$ be an arbitrary element. Then:*

1. $|hg|_A \geq |h|_A + |g|_A - K$;
2. *the path $[1, h] \cup h[1, g]$ is $K$-Hausdorff close to $[1, hg]$;*

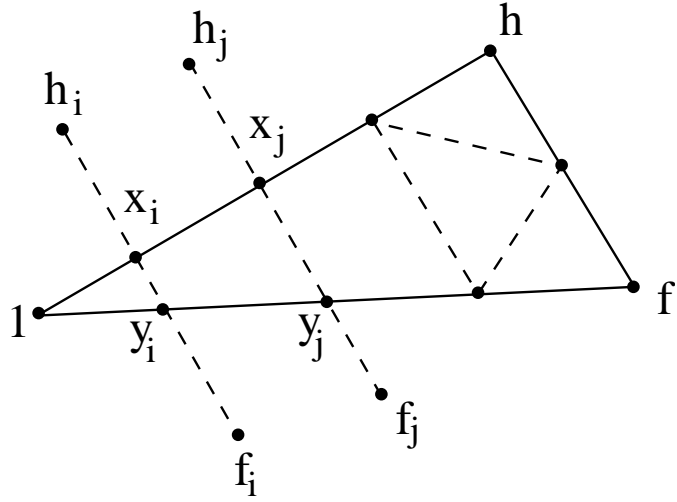

FIGURE 2. Figure for the proof of Proposition 4.5

3. $(g, h)_1 \leq K$.

Similarly, suppose $g \in G$ is shortest with respect to $d$ in the coset class $gH$. Let $h \in H$ be an arbitrary element. Then:

1. $|gh|_A \geq |g|_A + |h|_A - K$;
2. the path $[1, g] \cup g[1, h]$ is $K$-Hausdorff close to $[1, gh]$.
3. $(g^{-1}, h)_1 \leq K$.

**Proposition 4.5.** *Let $H, F \leq G$ be quasiconvex subgroups such that $H \cap F = 1$. Then $(H, F)_1 < \infty$.*

*Proof.* Let $C > 0$ be an integer such that both $H$ and $F$ are $C$-quasiconvex in $X$. Let $N$ be the number of elements $g \in G$ such that $|g|_A \leq 2C + \delta$.

Suppose $(H, F)_1 = \infty$. Then there are $h \in H, f \in F$ such that $(h, f)_1 \geq (N+1)(2C+1)$. Consider a geodesic triangle with vertices $1, h, f$ and sides $[1, h], [1, f]$ and $[h, f]$, as shown in Figure 2. Since $(h, f)_1 \geq (N+1)(2C+1)$, for any points $x \in [1, h], y \in [1, f]$ with $d(1, x) = d(1, y) \leq (N+1)(2C+1)$ we have $d(x, y) \leq \delta$. Define two sequences $x_i \in [1, h]$ and $y_i \in [1, f]$ for $i = 0, 1 \ldots, N$ so that $d(1, x_i) = d(1, y_i) = i(2C+1)$. Then $d(x_i, y_i) \leq \delta$ for each $i = 0, \ldots, N$. Since $H, F$ are $C$-quasiconvex, for each $i$ there are $h_i \in H, f_i \in F$ such that $d(x_i, h_i) \leq C$ and $d(y_i, f_i) \leq C$. Hence $d(h_i, f_i) = |f_i^{-1} h_i|_A \leq 2C + \delta$. Note that $d(h_i, h_j) \geq d(x_i, x_j) - 2C = |i-j|(2C+1) - 2C$. Hence for $i \neq j$ we have $d(h_i, h_j) > 0$ and $h_i \neq h_j$.

By the choice of $N$ there are some $i < j$ such that $e = f_i^{-1} h_i = f_j^{-1} h_j$. Hence $e(h_i^{-1} h_j)e^{-1} = f_i^{-1} f_j$. Therefore
$$f_i^{-1} h_i (h_i^{-1} h_j) h_i^{-1} f_i = f_i^{-1} f_j \Rightarrow h_j h_i^{-1} = f_j f_i^{-1} \neq 1,$$
contrary to our assumption that $H \cap F = 1$. □

**Definition 4.6** (Conjugacy separated subgroups)**.** Let $H, F \leq G$ be two subgroups. We say that $H$ and $F$ are *conjugacy separated in $G$* if for any $g \in G$
$$g^{-1}Hg \cap F = 1.$$

**Lemma 4.7.** *Let $F \leq G$ be a quasiconvex subgroup which is conjugacy separated from $H$ in $G$. There is an integer constant $K_1 = K_1(H, F, G, A) > 0$ with the following property.*

*Suppose $g \in G$ is shortest in the double coset $HgF$. Let $f \in F$ be an arbitrary element. Then for any $h \in H$ such that $hgf$ is shortest in $Hgf$ we have $|h|_A \leq K_1$.*

*Proof.* Let $K(H) > 0$ and $K(F) > 0$ be the integer constants provided by Lemma 4.4. Put $K := \max\{K(H), K(F)\}$. Let $C > 0$ be an integer such that both $H$ and $F$ are $C$-quasiconvex in $X$. Let $N$ be the number of elements in $G$ of length at most $2(K + C + \delta)$. Put $K_1 := N(2K + 2C + 2\delta + 1) + 3K + 2\delta$.

Suppose $f \in F$ and $h' \in H$ are such that $g' = hgf$ is shortest in $Hgf$. Thus $hg' = gf$, where $h = (h')^{-1}$.

Let $\alpha = [1, gf]$. Consider the geodesic triangles $\alpha \cup [1, g] \cup g[1, f]$ and $\alpha \cup [1, h] \cup h[1, g']$. By Lemma 4.4 there are points $p, q \in \alpha$ such that $d(p, g) \leq K$ and $d(h, q) \leq K$. Then $(g, p)_1 \geq d(1, p) - K$ and $d(h, q) \geq d(1, q) - K$. Let $a \in \alpha$ be such that $d(1, a) = \max\{d(1, p), d(1, q)\}$. By Lemma 3.7 the Gromov product is monotone non-decreasing and hence $(g, a)_1 \geq d(1, p) - K$, $(h, a)_1 \geq d(1, q) - K$. Hence by Proposition 3.6 and Lemma 4.4

$$K \geq (g, h)_1 \geq \min\{(g, a)_1, (h, a)_1\} - 2\delta \geq \min\{d(1, p), d(1, q)\} - K - 2\delta$$

and hence $\min\{d(1, p), d(1, q)\} \leq 2K + 2\delta$.

Suppose first that $d(1, q) \leq d(1, p)$. Then $|h|_A = d(1, h) \leq d(1, q) + K \leq 3K + 2\delta \leq K_1$, as required.

Suppose now that $d(1, p) \leq d(1, q)$, as shown in Figure 3. Hence $d(1, p) \leq 2K + 2\delta$ and therefore $|g|_A = d(1, g) \leq 3K + 2\delta$. Suppose $d(p, q) > N(2K + 2C + 2\delta + 1)$.

Choose the points $x_0 = p, x_1, \ldots, x_N$ on $[p, q] \subseteq \alpha$ so that $d(p, x_i) = i(2K + 2C + 2\delta + 1)$. Since $d(g, p) \leq K$, for each $i$ there is $y_i \in g[1, f]$ with $d(x_i, y_i) \leq K + \delta$. Since $F$ is $C$-quasiconvex, for each $i$ there is $f_i \in F$ such that $d(y_i, gf_i) \leq C$. Similarly, since $d(h, q) \leq K$, for each $i$ there is $z_i \in [1, h]$ with $d(x_i, z_i) \leq K + \delta$. Since $H$ is $C$-quasiconvex, for each $i$ there is $h_i \in H$ such that $d(z_i, h_i) \leq C$. Thus for every $i = 0, \ldots, N$ we have $|h_i^{-1}gf_i|_A = d(h_i, gf_i) \leq 2(C + K + \delta)$. Moreover for $i \neq j$ we have

$$d(h_i, h_j) \geq d(x_i, x_j) - 2C - 2K - 2\delta =$$
$$= |i - j|(2K + 2C + 2\delta + 1) - 2C - 2K - 2\delta > 0$$

and hence $h_i \neq h_j$.

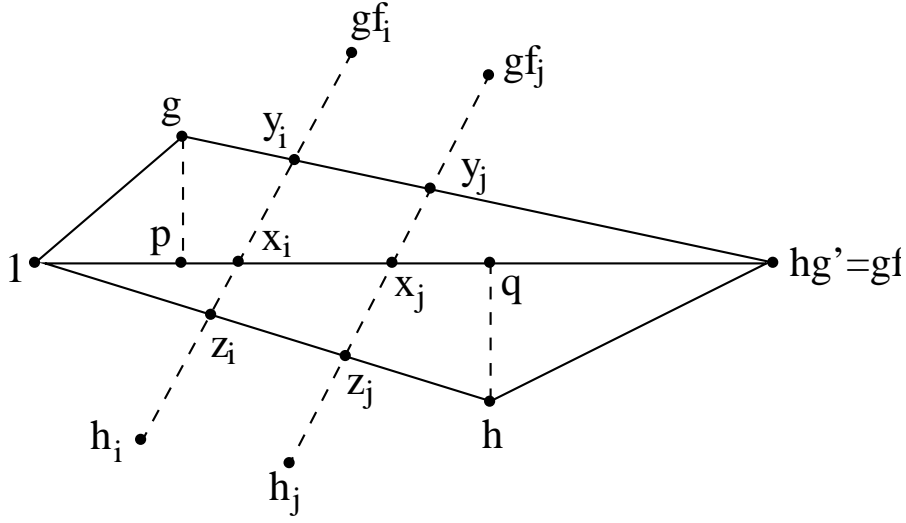

Figure 3. Figure for the proof of Lemma 4.7

By the choice of $N$ there are some $i < j$ such that
$$e = h_i^{-1} g f_i = h_j^{-1} g f_j \in G.$$

Hence
$$e(f_i^{-1} f_j) e^{-1} = h_i^{-1} h_j \in H - \{1\},$$

contrary to the assumption that $F$ and $H$ are conjugacy separated in $G$.

Thus $d(p, q) \leq N(2K + 2C + 2\delta + 1)$. Since $d(1, p) \leq 2K + 2\delta$ and $d(h, q) \leq K$, this implies
$$|h|_A = d(h, 1) \leq d(h, q) + d(h, p) + d(p, 1) \leq$$
$$\leq N(2K + 2C + 2\delta + 1) + 3K + 2\delta = K_1,$$

as required. □

**Proposition 4.8.** *Suppose $H, F \leq G$ are conjugacy-separated quasiconvex subgroups of $G$. Then there exists an integer constant $K_2 > 0$ with the following property.*

*Suppose $g \in G$ is shortest in the double coset $HgF$. Suppose $h \in H$, $f_1, f_2 \in F$ are such that $g f_1 = h g f_2$. Put $f = f_1 f_2^{-1} = g^{-1} h g$, so that $gf = hg$ and $f_1 = f f_2$. Then $|h|_A \leq K_2$ and $|f|_A = d(g, hg) \leq K_2$.*

*Proof.* Let $C > 0$ be an integer such that $H$ and $F$ are quasiconvex in $X = \Gamma(G, A)$. Let $K(H), K(F) > 0$ be the integer constants provided by Lemma 4.4. Put $K = \max\{K(H), K(F)\}$. Let $K_1 = K_1(H, F) > 0$ be the integer constant provided by Lemma 4.7.

For $i = 1, 2$ let $h_i \in F$ be such that $g_i' := h_i g f_i$ is shortest in $H g f_i$. Then by Lemma 4.7 we have $|h_1|_A, |h_2|_A \leq K_1$. From $g f_1 = h g f_2$ we have $h_1 g_1' = h h_2 g_2'$ and hence $g_1' = h_1^{-1} h h_2 g_2'$. By Lemma 4.4 we have

$|h_1^{-1}hh_2|_A \le K$. Hence $|h|_A \le K + 2K_1$. Fix a geodesic $\alpha = [1, gf_1] = [1, hgf_2]$. By Lemma 4.4 there is a point $x \in h[1, gf_2]$ such that $d(hg, x) \le K$. Since $|h|_A \le K + 2K_1$, there is a point $y \in \alpha$ such that $d(x, y) \le |h|_A + \delta \le K + 2K_1 + \delta$. Hence $d(hg, y) \le 2K + 2K_1 + \delta$. Since $|h|_A \le K + 2K_1$, we have

$$|d(1, y) - d(h, hg)| = |d(1, y) - |g|_A| \le 3K + 4K_1 + \delta.$$

On the other hand Lemma 4.4 implies that there is $z \in \alpha$ with $d(g, z) \le K$. Hence $|d(1, z) - |g|_A| \le K$. Since both $z, y$ are on $\alpha$, this implies

$$d(y, z) \le 4K + 4K_1 + \delta$$

Hence

$$d(hg, g) \le d(hg, y) + d(y, z) + d(z, g) \le$$
$$\le (2K + 2K_1 + \delta) + (4K + 4K_1 + \delta) + K = 7K + 6K_1 + 2\delta.$$

Thus the statement of Proposition 4.8 holds with $K_2 = 7K + 6K_1 + 2\delta$. □

Our main technical tool will be the following statement, whose proof is postponed until the last section.

**Theorem 4.9.** *Let $G$ be a non-elementary hyperbolic group and let $H \le G$ be a quasiconvex subgroup of infinite index. Then there exist elements $a, b \in G$ such that the elements $a, b$ generate a free quasiconvex subgroup of rank two which is conjugacy separated from $H$ in $G$.*

## 5. Proof of the main result

*Proof of Theorem 1.2.* If $H = 1$ then $G = G/H$ and the statement is obvious since $G$ contains a free subgroup of rank two and hence $G$ is non-amenable. Suppose $H \ne 1$. Let $X = \Gamma(G, A)$ be the Cayley graph of $G$ with respect to $A$. We will denote the word-metric on $X$ corresponding to $A$ by $d_X$. Also for $g \in G$ we denote $|g|_X := d_X(1, g)$. Put $Y = \Gamma(G, H, A)$. Thus $Y$ is a connected $2k$-regular graph where $k$ is the number of elements in $A$. We denote the simplicial metric on $Y$ by $d_Y$. By Theorem 4.9 there exists a quasiconvex free subgroup $F = F(a, b) \le G$ which is conjugacy separated from $H$ in $G$. Let $M = \max\{|a|_A, |b|_A\}$. Put $Q := \{a, b, a^{-1}, b^{-1}\}$. Recall that $H, F \le G$ are quasiconvex. Let $K(H) > 0, K(F) > 0$ be the integer constant provided by Lemma 4.4. Put $K = \max\{K(H), K(F)\}$. Let $K_2 = K_2(H, F) > 0$ be the integer constant provided by Proposition 4.8. Let $Z$ be the Cayley graph of $F$ with respect to the generating set $\{a, b\}$. Thus $Z$ is a 4-regular tree. Let $B$ be the set of all elements $f \in F$ such that $|f|_A \le K_2$. Let $T = |B|$ be the number of elements of $B$.

Let $S$ be an arbitrary finite nonempty subset of $\{Hg \mid g \in G\} = VY$. We will first decompose $S$ according to the double $H - F$-cosets of its members, that is write $S = S_1 \cup \cdots \cup S_n$, where $n \ge 1$ and

1. Each $S_i$ has the form $S_i = Hg_iF_i$, where $F_i \subseteq F$ is a finite nonempty subset and where $g_i \in G$ is shortest with respect to $d$ in the double coset $Hg_iF$;
2. For $i \neq j$ we have $Hg_iF \neq Hg_jF$.

We enumerate each $F_i$ as $F_i = \{f_{i,1}, \ldots, f_{i,m_i}\}$, where $m_i \geq 1$ is the number of elements in $F_i$.

**Claim.** Let $f' \in F$ be an arbitrary fixed element. Then for any $f'' \in F$ such that $Hg_if' = Hg_if''$ there is $f \in B$ such that $f'' = ff'$. Thus for any nonempty subset $R \subseteq F$ the set $\{Hg_if' \,|\, f' \in R\} \subseteq VY = \{Hg \,|\, g \in G\}$ consists of at least $|R|/T$ distinct elements.

Indeed, suppose that $f'' \in F$ is such that $Hg_if' = Hg_if''$. Lemma 4.8 implies that there is $h \in H$, $f \in F$ with $|h|_A \leq K_2, |f|_A \leq K_2$ such that $g_if'' = hg_if'$, $g_if = hg_i$ and $ff' = f''$. Since $|f|_A \leq K_2$, then $f \in B$ by definition of $B$. Also, by the choice of $T = |B|$ for a fixed element $f' \in F$ there are at most $T$ elements $f'' \in F$ which can arise in this fashion.

Since $Z$ is a 4-regular tree, it is non-amenable. Hence by part 4 Proposition 2.3 there is an integer $m > 0$ such that for any finite nonempty subset $R \subseteq VZ$ we have $|\mathcal{N}_m^Z(R)| \geq 2T|R|$.

Hence $|\mathcal{N}_m^Z(F_i)| \geq 2T|F_i| = 2T|S_i|$. Recall that $M = \max\{|a|_A, |b|_A\} = \max_{q \in Q} |q|_A$. Thus for any $f \in \mathcal{N}_m^Z(F_i)$ we have $Hg_if \in \mathcal{N}_{Mm}^Y(S_i)$. The Claim implies that $\{Hg_if \,|\, f \in \mathcal{N}_m^Z(F_i)\} \subseteq VY = \{Hg \,|\, g \in G\}$ consists of at least $\frac{1}{T}|\mathcal{N}_m^Z(F_i)| \geq \frac{1}{T}2T|S_i| = 2|S_i|$ distinct elements. Since $Hg_iF \cap Hg_jF = \emptyset$ for $i \neq j$, we have:

$$|\mathcal{N}_{mM}^Y(S)| \geq \sum_{i=1}^n 2|S_i| = 2\sum_{i=1}^n |S_i| = 2|S|.$$

Since $M$ and $m$ are fixed and a finite set $S \subseteq VY$ was chosen arbitrarily, the graph $Y$ is non-amenable by part 3 Proposition 2.3, and the theorem is proved. $\square$

We can now obtain Corollary 1.4 from the Introduction.

**Corollary 5.1.** *Let $G = \langle x_1, \ldots, x_k \,|\, r_1, \ldots, r_m \rangle$ be a non-elementary hyperbolic group and let $H \leq G$ be a quasiconvex subgroup of infinite index. Let $a_n$ be the number of freely reduced words in $A = \{x_1, \ldots, x_k\}^{\pm 1}$ of length $n$ representing elements of $H$. Let $b_n$ be the number of all words in $A$ of length $n$ representing elements of $H$.*

*Then*

$$\limsup_{n \to \infty} \sqrt[n]{a_n} < 2k - 1$$

*and*

$$\limsup_{n \to \infty} \sqrt[n]{b_n} < 2k.$$

*Proof.* Note that $k \geq 2$ since $G$ is non-elementary. Put $A = \{x_1, \ldots, x_k\}$ and $Y = \Gamma(G, H, A)$. We choose $x_0 := H1 \in VY$ as the base-vertex of $Y$. Note that $Y$ is 2k-regular by construction. Also, for any vertex $x$ of $Y$ and any word $w$ in $A \cup A^{-1}$ there is a unique path in $Y$ with label $w$ and origin $x$. The definition of Schreier subgroup graphs also implies that:

(1) A freely reduced word $w$ represents an element of $H$ if and only if the path in $Y$ labeled $w$ with origin $x_0$ terminates at $x_0$.

(2) A word $w$ represents an element of $H$ if and only if the path in $Y$ labeled $w$ with origin $x_0$ terminates at $x_0$.

Therefore $a_n(Y)$ equals the number of freely reduced words in the alphabet $A = \{x_1, \ldots, x_k\}^{\pm 1}$ of length $n$ representing elements of $H$. Similarly, $b_n(Y)$ equals the number of all words in $A$ of length $n$ representing elements of $H$.

By Theorem 1.2 $Y$ is non-amenable. Hence by Theorem 2.5 $\alpha(Y) < 2k-1$ and $\beta(Y) < 2k$, as required. $\square$

## 6. Constructing a conjugacy separated cyclic subgroup

We will now concentrate on proving Theorem 4.9.

Until the end of this article, unless specified otherwise, let $G$ be a non-elementary hyperbolic group and let $H \leq G$ be an infinite quasiconvex subgroup of infinite index. Fix a finite generating set $A$ of $G$ defining the word metric $d$ on $X = \Gamma(G, A)$. Let $\delta \geq 1$ be an integer such that all geodesic triangles in $X$ are $\delta$-trim.

The following lemma is a straightforward hyperbolic exercise and we leave the proof to the reader:

**Lemma 6.1.** *The following hold in $X$:*

1. *Let $x, y \in X$ and suppose that $p \in [1, x]$, $q \in [1, y]$ are such that $d(1, p), d(1, q) \geq (p, q)_1 + 10\delta$. Then $(x, y)_1 \leq (p, q)_1 + 2\delta$.*
2. *Suppose $x, g_1, g_2 \in G$ and $(g_1, xg_2)_1 = L \leq |g_2|_A$. Then for any points $y \in [1, g_1]$ and $z \in x[1, g_2]$ with $d(1, y) = d(1, z) \leq L$ we have $d(y, z) \leq |x|_A + 2\delta$.*
3. *Let $p, x, y, q \in X$ be such that $T_1 = (y, p)_x$, $T_2 = (x, q)_y$ and that $d(x, y) \geq T_1 + T_2 + 10\delta$. Then $d(p, q) \geq d(p, x) + d(q, y) - T_1 - T_2$.*
4. *Suppose $x, y, z \in X$ and $(x, z)_y = T$. Then any path $[x, y] \cup [y, z]$ is $(1, L)$-quasigeodesic for some constant $L = L(T, \delta)$.*
5. *For any $L_1, L_2 > 0$ there is a constant $L_3 > 0$ with the following property. Suppose $x \in G$ with $|x|_A \leq L_1$. suppose $g_1, g_2 \in G$ are such that $(g_1, xg_2)_1 \leq L_2$. Then the paths $[g_1^{-1}, 1] \cup [1, xg_2]$ and $[g_1^{-1}, 1] \cup [1, x] \cup x[1, g_2]$ are $(1, L_3)$-quasigeodesics.*

**Definition 6.2.** Let $\lambda \geq 1$. We say that an element $g \in G$ is $\lambda$-*cyclically reduced* if there exists a word $w$ representing $g$ such that every cyclic permutation of $w$ is $\lambda$-quasigeodesic. We will say that $g \in G, g \neq 1$ is *periodically geodesic* if for any geodesic representative $w$ of $g$ and any integer $n > 0$ the word $w^n$ is also geodesic.

The facts that the language of geodesics on $G$ is regular and every cyclic subgroup is quasiconvex easily imply:

**Lemma 6.3.** [4] *Let $g \in G$ be an element of infinite order. Then for some $n > 0$ the element $g^n$ is conjugate to a periodically-geodesic element.*

**Proposition 6.4.** *For any $x \in G$ of infinite order and any $L > 0$ there exist a constant $T = T(x, L)$ such that $\lim_{L \to \infty} T(x, L) = \infty$ and with the following property.*

*Suppose $x \in G$ and suppose $g \in G$ is such that $(xg, g)_1 \geq L$. Then*

$$(g, \langle x \rangle)_1 \geq T.$$

*Proof.* Let $\lambda > 0$ be such that for any $n$ we have $|x^n| \geq \lambda |n|$. Put

$$T := \frac{2L - \frac{(10\delta + \lambda)(|x|_A + 2\delta)}{\lambda}}{2 + \frac{2|x|_A + 4\delta}{\lambda}} - 1.$$

Suppose $(g, \langle x \rangle)_1 < T$.
Let

$$m := \lfloor \frac{2T + 10\delta + \lambda}{\lambda} \rfloor.$$

Then $|x^m|_A \geq \lambda m \geq 2T + 10\delta$. The assumption $(g, \langle x \rangle)_1 < T$ implies that $(g, x^m)_1 \leq T$ and $(1, x^m g)_{x^m} < T$. Hence by part 3 of Lemma 6.1 for any $p \in [1, g]$ and $q \in x^m[1, g]$ we have $d(p, q) \geq d(p, 1) + d(q, x^m) - 2T$.

Let $p \in [1, g]$ be such that $d(1, p) = L$.

Then by part 2 of Lemma 6.1 $d(p, xp) \leq |x|_A + 2\delta$ and hence inductively $d(p, x^m p) \leq m|x|_A + 2m\delta$. By the previous remark $d(p, x^m p) \geq 2L - 2T$.

Therefore

$$2L - 2T \leq d(p, x^m p) \leq m|x|_A + 2m\delta = m(|x|_A + 2\delta) \Rightarrow$$
$$2L - 2T \leq \frac{2T + 10\delta + \lambda}{\lambda}(|x|_A + 2\delta) \Rightarrow$$
$$2L \leq T(2 + \frac{2|x|_A + 4\delta}{\lambda}) + \frac{(10\delta + \lambda)(|x|_A + 2\delta)}{\lambda} \Rightarrow$$
$$T \geq \frac{2L - \frac{(10\delta + \lambda)(|x|_A + 2\delta)}{\lambda}}{2 + \frac{2|x|_A + 4\delta}{\lambda}}$$

which contradicts the definition of $T$. □

**Proposition 6.5.** *There exist integer constants $c' = c'(G, A) > 0$ and $L' = L'(G, A) > 0$ such that for any $g \in G$ there is an $L'$-cyclically reduced element $g'$ such that $d(g, g') \leq c'$.*

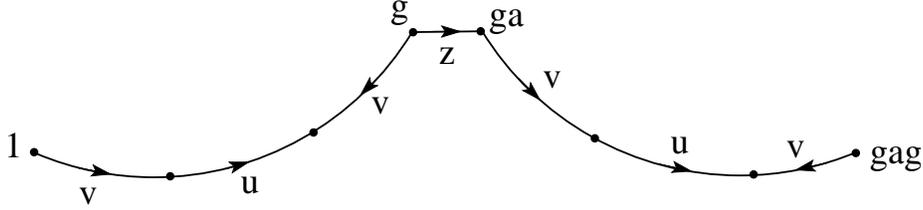

FIGURE 4. The word $wzw$ representing $gag$

*Proof.* Since $G$ is non-elementary and contains a free subgroup of rank two, there exist $a, b \in G$ such that $\langle a \rangle \cap \langle b \rangle = 1$. Since the cyclic subgroups $\langle a \rangle, \langle b \rangle \leq G$ are quasiconvex, Lemma 4.5 implies that $(\langle a \rangle, \langle b \rangle)_1 < \infty$.

Let $T_0 := (\langle a \rangle, \langle b \rangle)_1 + 2\delta + 1$. There exists an integer $L > 1$ such that $T(a, L) > T_0$ and $T(b, L) > T_0$, where $T(a, L), T(b, L)$ are provided by Proposition 6.4.

By part 5 of Lemma 6.1 there is an integer constant $L_0 > L + 100\delta > 1$ such that if $(h, xh)_1 \leq L$, where $x \in \{a, b\}$ and $h \in G$, then any path $[h^{-1}, 1] \cup [1, x] \cup x[1, h]$ is $(1, L_0)$-quasigeodesic. By Proposition 3.9 there are integers $N > L_0$ and $L_1 > L_0$ such that any $N$-local $(1, L_0)$-quasigeodesic is $L_1$-quasigeodesic.

Let $g \in G$ be an arbitrary element.

If $(g, g^{-1})_1 \leq 100L_0 + 100N + 100\delta$, then $[1, g] \cup g[1, g]$ is an $L_2$-quasigeodesic for some integer constant $L_2 > 1$ provided by part 3 of Lemma 6.1. This means that for any geodesic word $w$ representing $g$ the word $ww$ is $L_2$-quasigeodesic. Hence any cyclic permutation of $w$ is $L_2$-quasigeodesic and so $g$ is an $L_2$-cyclically reduced element.

Assume now that $(g, g^{-1})_1 > 100L_0 + 100N + 100\delta$. Then $g$ has a $(1, 100\delta)$-quasigeodesic representative $w$ of the form $w = vuv^{-1}$ where the words $u, v$ are geodesics and $|v| = 100L + 100N + 100\delta$. Put $f = \overline{v} \in G$. Recall that $L_0 > 100\delta$ by construction.

Suppose that $(f, af)_1 > L$ and $(f, bf) > L$. This implies that $(f, \langle a \rangle)_1 > T_0$ and $(f, \langle b \rangle)_1 > T_0$ by Lemma 6.4 and the choice of $L$. Therefore

$$(\langle a \rangle, \langle b \rangle)_1 \geq T_0 - 2\delta = (\langle a \rangle, \langle b \rangle)_1 + 1,$$

which is impossible. Thus either $(f, af)_1 \leq L$ or $(f, bf) \leq L$. Without loss of generality we may assume that $(f, af)_1 \leq L$. Let $z$ be a geodesic word representing $a$. Then by the choice of $L_0$ the word $v^{-1}zv$ is a $(1, L_0)$-quasigeodesic. Therefore the word $vuv^{-1}zvuv^{-1} = wzw$ is an $N$-local $(1, L_0)$-quasigeodesic and thus $L_1$-quasigeodesic by the choice of $L_1$ (see Figure 4). Hence any cyclic permutation of $wz$ is $L_1$-quasigeodesic and so $ga$ is an $L_1$-cyclically reduced element. Note that $|a|_X \leq c$.

Thus we have verified that the statement of Proposition 6.5 holds with $c' = \max\{|a|_A, |b|_A\}$ and $L' = \max\{L_1, L_2\}$. □

**Corollary 6.6.** *There are integer constants $c = c(G, A) > 0$ and $\lambda = \lambda(G, A) > 0$ with the following property. For any $g \in G$ there is an element of infinite order $a_g \in G$ such that $d(g, a_g) \le c$ and that for some word $w_g$ representing $a_g$ the word $w_g^n$ is $\lambda$-quasigeodesic for any integer $n > 0$.*

*Proof.* Let $c' > 0, L' > 0$ be the integer constants provided by Proposition 6.5. By Proposition 3.9 there exist integers $N > 1, \lambda > 1$ such that any $N$-local $L'$-quasigeodesic in $X$ is $\lambda$-quasigeodesic. Let $L', c' > 0$ be integer constants provided by Proposition 6.5.

Let $B$ be the set of all elements of $G$ of length at most $L'(N + L') + c'$. Choose $c'' > 0$ such that every element $g$ of $B$ is at most $c''$-away from a periodically geodesic element $a(g)$ of $G$.

Suppose now $g \in G$ and $g \notin B$, so that $|g|_A \ge L'(N + L') + c' + 1$. By proposition 6.5 there is an element $a(g) \in G$ with $d(a(g), g) \le c'$ such that $a(g)$ is $L'$-cyclically reduced. Hence $|a(g)|_X \ge L'(N + L')$. Moreover, there is an $L'$-quasigeodesic representative $w_g$ of $a(g)$ such that every cyclic permutation of $w_g$ is $L'$-quasigeodesic. Thus $|w_g| \ge N$. By the choice of $N$ this implies that for any $n > 1$ the word $w_g^n$ is $\lambda$-quasigeodesic. In particular $a(g)$ has infinite order in $G$.

Thus $\lambda$ and $c = \max\{c', c''\}$ satisfy the requirements of the lemma. □

The following useful lemma is due to B.H.Neumann [54]:

**Lemma 6.7.** *Let $G$ be a group. Suppose $G = a_1 H_1 \cup \ldots a_k H_k$ where $a_i \in G$ and $H_i \le G$ are subgroups of $G$.*
*Then at least one of $H_i$ has finite index in $G$.*

**Proposition 6.8.** *Let $G$ be a non-elementary hyperbolic group and let $H \le G$ be a quasiconvex subgroup of infinite index. Then there is an element $a \in G$ of infinite order such that the subgroup $\langle a \rangle$ is conjugacy separated from $H$ in $G$.*

*Proof.* Let $A$ be a fixed finite generating set of $G$ and let $X = \Gamma(G, A)$ be the Cayley graph of $G$ with the word-metric $d$. Let $\delta > 0$ be an integer such that all geodesic triangles in $X$ are $\delta$-trim. Let $C > 0$ be and integer such that $H$ is $C$-quasiconvex in $X$. Let $c > 0, \lambda > 0$ be the constants provided by Corollary 6.6.

Suppose that no infinite cyclic subgroup of $G$ is conjugacy separated from $H$ in $G$.

Let $g \in G$ be an arbitrary element. By Corollary 6.6 there exists an element of infinite order $a_g \in G$ with a representative word $w_g$ such that $d(g, a_g) \le c$ and such that the word $w_g^n$ is $\lambda$-quasigeodesic for any $n \ge 1$. Recall that $\lambda \ge 1$. By assumption there is some $n > 0$ such that $a_g^n = xhx^{-1}$ for some $x \in G$ and $h \in H$. Let $E > 0$ be such that any two $\lambda$-quasigeodesics with common endpoints in $X$ are $E$-close.

Choose $m > 1$ such that $|a_g^{mn}|_A > \lambda[2|x|_A + 2E + 2|w_g| + \lambda]$. We have $a_g^{mn} = xh^m x^{-1}$. By the choice of $m$ there is a subsegment $J = [p, p']$ of

$[1, a_g^{mn}]$ of length at least $\lambda[2E + 2|w_g| + \lambda]$ which is contained in the $2\delta$-neighborhood of $x[1, h]$.

Let $q, q'$ be the points on the path $W$ from 1 to $a_g^{mn}$ labeled by $w_g^{mn}$ such that $d(p, q) \leq E$ and $d(p', q') \leq E$. Hence $d(p', q') \geq \lambda[3|w_g| + \lambda]$. Therefore the length of the segment $U$ of $W$ from $q$ to $q'$ is at least $2|w_g|$. Hence $U$ contains some two vertices the form $v_0 = a_g^t, v_1 = a_g^{t+1}$. For every $i = 0, 1$ there is a point $s_i \in J$ with $d(v_i, s_i) \leq E + 2\delta$. By the choice of $J$ for every $i$ there is $z_i \in x[1, h]$ such that $d(s_i, z_i) \leq 2\delta$. Finally for each $z_i$ there is $h_i \in H$ such that $d(z_i, xh_i) \leq C$.

Thus $d(v_i, xh_i) = d(a_g^{t+i}, xh_i) \leq C + E + 4\delta$.

Put $f_0 = a_g^{-t} xh_0$ and $f_1 = h_1^{-1} x^{-1} a_g^{t+1}$. Then $|f_0|_A, |f_1|_A \leq C + E + 4\delta$ and $a_g = f_0 h_0^{-1} h_1 f_1$.

Recall that $d(a_g, g) \leq c$, so that $g = f_0 h_0^{-1} h_1 f_1'$ where $f_1' \leq C + E + 4\delta + c$. Since $g \in G$ was chosen arbitrarily, we have established that $G$ can be represented as a finite union

$$G = \cup \{aHb \,|\, a, b \in G, |a|_A, |b|_A \leq C + E + 4\delta + c\}$$

However, each set $aHb$ can be written as a coset $(ab)b^{-1}Hb$ of the subgroup $b^{-1}Hb \leq G$. Hence $G$ is the finite union of cosets of subgroups of infinite index, which is impossible by Lemma 6.7. □

## 7. Constructing a conjugacy separated free subgroup

Let $c \in G$ be an element of infinite order provided by Proposition 6.8, so that the infinite cyclic subgroup $L = \langle c \rangle$ is conjugacy separated from $H$ in $G$.

**Proposition 7.1.** *There exist constants $M' > 1, C' > 0, \lambda' > 0$ with the following properties.*

*Suppose $h_1, \ldots, h_k \in H$, $n_0, \ldots, n_k \in \mathbb{Z}$ are such that $k \geq 0$, such that $|h_i|_A \geq C'$ for $i = 1, \ldots, k$ and that $|n_i| \geq M'$ for $0 < i < k$. Let $u_i$ be an $A$-geodesic word representing $h_i$ and let $v_i$ be an $A$-geodesic word representing $c^{n_i}$. Then the word*

$$W = v_0 u_1 v_1 u_2 \ldots \ldots u_k v_k$$

*is a $\lambda'$-quasigeodesic in $X$.*

*Proof.* Since $L$ and $H$ are conjugacy separated in $G$, we have $L \cap H = 1$ and hence $(H, L)_1 < \infty$ by Lemma 4.5. Let $T_1 = (H, L)_1$. By part 3 of Lemma 6.1 there is an integer constant $T_2 > T_1$ such that whenever $x, y, z \in X$ are such that $(x, z)_y \leq T_1$ then $[x, y] \cup [y, z]$ is a $T_2$-quasigeodesic. By Proposition 3.9 there are integers $N > 0, \lambda' > 0$ such that any $N$-local $T_2$-quasigeodesic is $\lambda'$-quasigeodesic in $X$. Choose $M' > 0$ such that for any $n \geq M'$ we have $|c^n| \geq 100(N + T_2)$. Put $C' = 100(N + T_2)$. Then any word $W$ satisfying the conditions of Proposition 7.1 with these values of $C', M'$ is an $N$-local $T_2$-quasigeodesic and hence $\lambda'$-quasigeodesic, as required. □

We are now ready to prove Theorem 4.9.

*Proof of Theorem 4.9.* Fix a finite generating set $A$ of $G$ defining the word metric $d$ on $X = \Gamma(G, A)$. Let $\delta \geq 1$ be an integer such that all geodesic triangles in $X$ are $\delta$-trim. If $H$ is finite, the statement of the theorem is obvious since $G$ is non-elementary and hence contains a free quasiconvex subgroup of rank two. Thus we will assume that $H$ is infinite.

Let $c \in G$ be an element of infinite order provided by Proposition 6.8, so that the infinite cyclic subgroup $L = \langle c \rangle$ is conjugacy separated from $H$ in $G$. Fix $C > 0$ such that $H$ and $L$ are $C$-quasiconvex in $X$. Let $M', C', \lambda'$ be the constants provided by Proposition 7.1. Since $H \leq G$ is infinite, there is an element $h' \in H$ of infinite order so that $|(h')^n|_A \geq C'$ for any $n \neq 0$. Let $C'' > 0$ be such that any two $\lambda$-quasigeodesics in $X$ with common endpoints are $C''$-Hausdorff close. Let $N$ be the number of elements in $G$ of length at most $C'' + 2C + 2\delta$.

Choose $m \geq M$ such that for any integer $t$ with $|t| \, m$ we have $|c^t|_A > \lambda'^2 + 1$ and $|c^t| \geq 10(N+1)(2C + 4\delta + 1) + 10|h'|_A$. Let $C_1 > 0$ be the quasiconvexity constant of the subgroup $\langle c^m \rangle$ in $X = \Gamma(G, A)$.

Put $a := c^m, b := (h')^{-1} c^m h'$. We claim that $F := \langle a, b \rangle \leq G$ is a free quasiconvex subgroup of rank two. If $f$ is a nontrivial freely reduced product in $a, b$ then by Proposition 7.1 and the choice of $m$ we have $|f|_A \geq \frac{\lambda'^2+1}{/}\lambda' - \lambda' > 0$. Hence $f \neq 1$ and so $F$ is indeed free on $a, b$. Moreover, Proposition 7.1 also implies that for any $f \in F$ and any $p \in [1, f]$ there is $f' \in F$ such that $d(p, f') \leq C'' + C_1 + |h'|_A$. Thus $F$ is quasiconvex in $G$, as required.

We will now show that $F$ is conjugacy separated from $H$ in $G$. Indeed, suppose not. Then for some $g \in G$, $f \in F - \{1\}$ and $h \in H - \{1\}$ we have $ghg^{-1} = f$.

For any integer $n$ we have $gh^n g^{-1} = f^n$. Choose $n > 1$ so that $|h^n|_A > 2|g|_A + 2(N+1)(2C + 4\delta + 1) + 2|h'|_A + 4\delta$. Consider a geodesic quadrilateral in $X$ with vertices $1, g, gh^n, gh^n g^{-1} = f^n$ and sides $\alpha = [1, f^n]$, $[1, g]$, $\beta = g[1, h^n]$ and $gh^n[1, g^{-1}]$.

By the choice of $n$ there is a subsegment $J_1$ of $\beta$ of length at least $2(N+1)(2C + 4\delta + 1) + 2|h'|_A + 4\delta$ such that $J_1$ is $2\delta$-Hausdorff close to a subsegment $J_2$ of $\alpha$. Thus the length of $J_2$ is at least $2(N+1)(2C + 4\delta + 1) + 2|h'|_A$.

By the choice of $m$ Proposition 7.1 implies that there is a subsegment $J_2'$ of $J_2$ of length at least $(N+1)(2C + 4\delta + 1)$ such that $J_2'$ is contained in the $C''$-neighborhood of the path of the form $g_1[1, c^t]$ for some $g_1 \in G$ and some $t$ with $|t| \geq m$.

Consider a sequence of points $x_0, x_1, \ldots, x_N$ in $J_2'$ such that $d(x_i, x_j) = 2C + 4\delta + 1$, as shown in Figure 5.

For each $i = 0, \ldots, N$ there is a point $y_i \in J_1 \subseteq \beta$ such that $d(y_i, x_i) \leq 2\delta$. Since $H$ is $C$-quasiconvex, for each $i$ there is $h_i \in H$ such that $d(y_i, gh_i) \leq C$. Note that for $i \neq j$ we have

$$d(gh_i, gh_j) \geq d(x_i, x_j) - 2C - 4\delta = |i - j|(2C + 4\delta + 1) - 2C - 4\delta > 0$$

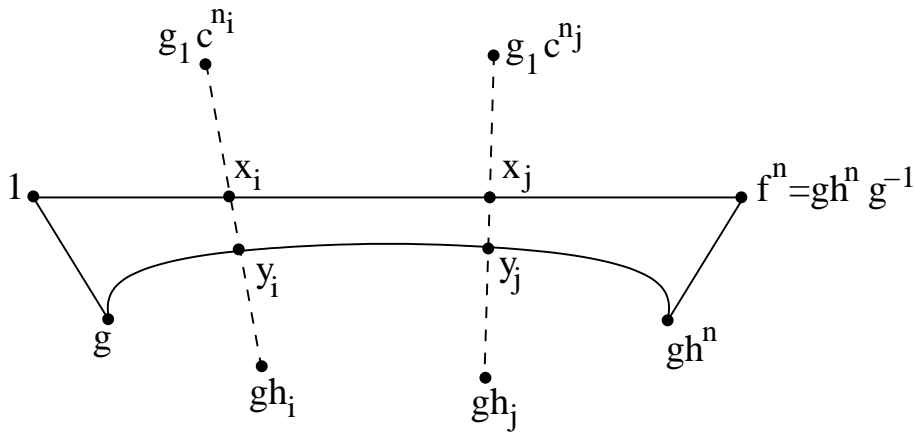

Figure 5. Figure for the proof of Theorem 4.9

by the choice of $x_0, \ldots, x_N$. Hence $h_i \neq h_j$ for $i \neq j$.

Also for every $i$ there is $n_i$ such that $d(x_i, g_1 c^{n_i}) \leq C'' + C$, since the subgroup $L = \langle c \rangle$ is $C$-quasiconvex in $X$. Hence $d(g_1 c^{n_i}, gh_i) \leq C'' + 2C + 2\delta$. Therefore by the choice of $N$ there are some $i < j$ such that $e = h_i^{-1} g^{-1} g_1 c^{n_i} = h_j^{-1} g^{-1} g_1 c^{n_j} \in G$. Hence $ec^{n_j - n_i} e^{-1} = h_i^{-1} h_j \in H - \{1\}$, which contradicts the fact that $H$ and $\langle c \rangle$ are conjugacy separated in $G$.

This completes the proof of Theorem 4.9. □

## References

[1] J.Alonso, T.Brady, D.Cooper, V.Ferlini, M.Lustig, M.Mihalik, M.Shapiro and H.Short, *Notes on hyperbolic groups,* In: " Group theory from a geometrical viewpoint", Proceedings of the workshop held in Trieste, É. Ghys, A. Haefliger and A. Verjovsky (editors). World Scientific Publishing Co., 1991

[2] A. Ancona, *Théorie du potentiel sur les graphes et les variétés,* Ecole d'été de Probabilités de Saint-Flour XVIII—1988, 1–112, Lecture Notes in Math., **1427**, Springer, Berlin, 1990

[3] G. Arzhantseva, *On Quasiconvex Subgroups of Word Hyperbolic Groups,* Geometriae Dedicata **87** (2001), 191–208

[4] G. Baumslag, S. Gersten, M. Shapiro and H. Short, *Automatic groups and amalgams,* J. of Pure and Appl. Algebra **76** (1991), 229–316

[5] L. Bartholdi, *Counting paths in graphs.* Enseign. Math. (2) **45** (1999), no. 1-2, 83–131.

[6] I. Benjamini, R. Lyons and O. Schramm, *Percolation perturbations in potential theory and random walks,* Random walks and discrete potential theory (Cortona, 1997), 56–84, Sympos. Math., XXXIX, Cambridge Univ. Press, Cambridge, 1999

[7] I. Benjamini and O. Schramm, *Every graph with a positive Cheeger constant contains a tree with a positive Cheeger constant,* Geom. Funct. Anal. **7** (1997), no. 3, 403–419

[8] M. Bestvina and M. Feighn, *A combination theorem for negatively curved groups,* J. Differential Geom. **35** (1992), no. 1, 85–101

[9] B. Bowditch, *Cut points and canonical splittings of hyperbolic groups,* Acta Math. **180** (1998), no. 2, 145–186

[10] P. Bowers, *Negatively curved graph and planar metrics with applications to type,* Michigan Math. J. **45** (1998), no. 1, 31–53

[11] P. Brinkmann, *Hyperbolic automorphisms of free groups,* Geom. Funct. Anal. **10** (2000), no. 5, 1071–1089

[12] J. Block and S. Weinberger, *Aperiodic tilings, positive scalar curvature and amenability of spaces,* J. Amer. Math. Soc. **5** (1992), no. 4, 907–918

[13] J. W. Cannon, *The theory of negatively curved spaces and groups.* Ergodic theory, symbolic dynamics, and hyperbolic spaces (Trieste, 1989), 315–369, Oxford Sci. Publ., Oxford Univ. Press, New York, 1991

[14] J. Cao, *Cheeger isoperimetric constants of Gromov-hyperbolic spaces with quasi-poles,* Commun. Contemp. Math. **2** (2000), no. 4, 511–533

[15] T. Ceccherini-Silberstein, R. Grigorchuck and P. de la Harpe, *Amenability and paradoxical decompositions for pseudogroups and discrete metric spaces,* (Russian) Tr. Mat. Inst. Steklova **224** (1999), Algebra. Topol. Differ. Uravn. i ikh Prilozh., 68–111; translation in Proc. Steklov Inst. Math., **224** (1999), no. 1, 57–97

[16] F. R. K. Chung, *Laplacians of graphs and Cheeger's inequalities.* Combinatorics, Paul Erdös is eighty, Vol. 2 (Keszthely, 1993), 157–172, Bolyai Soc. Math. Stud., **2**, János Bolyai Math. Soc., Budapest, 1996

[17] F. R. K. Chung and K. Oden, *Weighted graph Laplacians and isoperimetric inequalities,* Pacific J. Math. **192** (2000), no. 2, 257–273

[18] J. Cohen, *Cogrowth and amenability of discrete groups.* J. Funct. Anal. **48** (1982), no. 3, 301–309

[19] M. Coornaert, T. Delzant, and A. Papadopoulos, *Géométrie et théorie des groupes. Les groupes hyperboliques de Gromov.* Lecture Notes in Mathematics, 1441; Springer-Verlag, Berlin, 1990

[20] J. Dodziuk, *Difference equations, isoperimetric inequality and transience of certain random walks,* Trans. Amer. Math. Soc. **284** (1984), no. 2, 787–794

[21] M. Dunwoody and MSageev, *JSJ-splittings for finitely presented groups over slender groups,* Invent. Math. **135** (1999), no. 1, 25–44

[22] M. Dunwoody and E. Swenson, *The algebraic torus theorem,* Invent. Math. **140** (2000), no. 3, 605–637

[23] G. Elek, *Amenability, $l_p$-homologies and translation invariant functionals,* J. Austral. Math. Soc. Ser. A **65** (1998), no. 1, 111–119

[24] D. Epstein, J. Cannon, D. Holt, S. Levy, M. Paterson, W. Thurston, *Word Processing in Groups,* Jones and Bartlett, Boston, 1992

[25] D. Epstein, and D. Holt, *Efficient computation in word-hyperbolic groups.* Computational and geometric aspects of modern algebra (Edinburgh, 1998), 66–77, London Math. Soc. Lecture Note Ser., **275**, Cambridge Univ. Press, Cambridge, 2000

[26] R. Foord, *Automaticity and Growth in Certain Classes of Groups and Monoids,* PhD Thesis, Warwick University, 2000

[27] K. Fujiwara and P. Papasoglu, *JSJ decompositions and complexes of groups,* preprint, 1996

[28] V. N. Gerasimov, *Semi-splittings of groups and actions on cubings,* in "Algebra, geometry, analysis and mathematical physics (Novosibirsk, 1996)", 91–109, 190, Izdat. Ross. Akad. Nauk Sib. Otd. Inst. Mat., Novosibirsk, 1997

[29] P. Gerl, *Amenable groups and amenable graphs.* Harmonic analysis (Luxembourg, 1987), 181–190, Lecture Notes in Math., **1359**, Springer, Berlin, 1988

[30] S. Gersten and H. Short, *Rational subgroups of biautomatic groups,* Ann. Math. (2) **134** (1991), no. 1, 125–158

[31] E. Ghys and P. de la Harpe (editors), *Sur les groupes hyperboliques d'aprés Mikhael Gromov,* Birkhäuser, Progress in Mathematics series, vol. **83**, 1990.

[32] R. Gitik, *On the combination theorem for negatively curved groups,* Internat. J. Algebra Comput. **6** (1996), no. 6, 751–760
[33] R. Gitik, *On quasiconvex subgroups of negatively curved groups,* J. Pure Appl. Algebra **119** (1997), no. 2, 155–169
[34] R. Gitik, *On the profinite topology on negatively curved groups,* J. Algebra **219** (1999), no. 1, 80–86
[35] R. Gitik, *Doubles of groups and hyperbolic LERF 3-manifolds,* Ann. of Math. (2) **150** (1999), no. 3, 775–806
[36] R. Gitik, *Tameness and geodesic cores of subgroups,* J. Austral. Math. Soc. Ser. A **69** (2000), no. 2, 153–16
[37] R. Gitik, M. Mitra, E. Rips, M. Sageev, *Widths of subgroups,* Trans. Amer. Math. Soc. **350** (1998), no. 1, 321–329
[38] R. I. Grigorchuk, *Symmetrical random walks on discrete groups.* Multicomponent random systems, pp. 285–325, Adv. Probab. Related Topics, **6**, Dekker, New York, 198
[39] M. Gromov, *Hyperbolic Groups*, in "Essays in Group Theory (G.M.Gersten, editor)", MSRI publ. **8**, 1987, 75–263
[40] M. Gromov, *Asymptotic invariants of infinite groups.* Geometric group theory, Vol. 2 (Sussex, 1991), 1–295, London Math. Soc. Lecture Note Ser., **182**, Cambridge Univ. Press, Cambridge, 1993
[41] D. Holt, *Automatic groups, subgroups and cosets.* The Epstein birthday schrift, 249–260, Geom. Topol. Monogr., **1**, Geom. Topol., Coventry, 1998
[42] V. Kaimanovich, *Equivalence relations with amenable leaves need not be amenable.* Topology, ergodic theory, real algebraic geometry, 151–166, Amer. Math. Soc. Transl. Ser. 2, **202**, Amer. Math. Soc., Providence, RI, 2001
[43] V. Kaimanovich and W. Woess, *The Dirichlet problem at infinity for random walks on graphs with a strong isoperimetric inequality.* Probab. Theory Related Fields **91** (1992), no. 3-4, 445–466
[44] I. Kapovich, *Detecting quasiconvexity: algorithmic aspects.* Geometric and computational perspectives on infinite groups (Minneapolis, MN and New Brunswick, NJ, 1994), 91–99; DIMACS Ser. Discrete Math. Theoret. Comput. Sci., **25**, Amer. Math. Soc., Providence, RI, 1996
[45] I. Kapovich, *Quasiconvexity and amalgams,* Internat. J. Algebra Comput. **7** (1997), no. 6, 771–811
[46] I. Kapovich, *A non-quasiconvexity embedding theorem for word-hyperbolic groups,* Math. Proc. Cambridge Phil. Soc. **127** (1999), no. 3, 461–486
[47] I. Kapovich, *The combination theorem and quasiconvexity,* Internat. J. Algebra Comput. **11** (2001), no. 2, 185–216.
[48] I. Kapovich, *The geometry of relative Cayley graphs for subgroups of hyperbolic groups,* preprint, 2002
[49] I. Kapovich, and H. Short, *Greenberg's theorem for quasiconvex subgroups of word hyperbolic groups,* Canad. J. Math. **48** (1996), no. 6, 1224–1244
[50] A. Lubotzky, *Discrete groups, expanding graphs and invariant measures.* With an appendix by Jonathan D. Rogawski. Progress in Mathematics, **125**, Birkhäuser Verlag, Basel, 1994
[51] M. Mihalik, *Group extensions and tame pairs,* Trans. Amer. Math. Soc. **351** (1999), no. 3, 1095–1107
[52] M. Mihalik and W. Towle, *Quasiconvex subgroups of negatively curved groups,* J. Pure Appl. Algebra **95** (1994), no. 3, 297–301
[53] M. Mitra, *Cannon-Thurston maps for trees of hyperbolic metric spaces,* J. Differential Geom. **48** (1998), no. 1, 135–164
[54] B. H. Neumann, *Groups covered by finitely many cosets,* Publ. Math. Debrecen **3** (1954), 227–242

[55] L. Reeves, *Rational subgroups of cubed* 3-*manifold groups,* Michigan Math. J. **42** (1995), no. 1, 109–126

[56] E. Rips, *Subgroups of small cancellation groups,* Bull. London Math. Soc. **14** (1982), no. 1, 45–47

[57] E. Rips and Z. Sela, *Cyclic splittings of finitely presented groups and the canonical JSJ decomposition,* Ann. of Math. (2) **146** (1997), no. 1, 53–109

[58] M. Sageev, *Ends of group pairs and non-positively curved cube complexes,* Proc. London Math. Soc. (3) **71** (1995), no. 3, 585–617

[59] M. Sageev, *Codimension-*1 *subgroups and splittings of groups,* J. Algebra **189** (1997), no. 2, 377–389

[60] R. Schonmann, *Multiplicity of phase transitions and mean-field criticality on highly non-amenable graphs,* Comm. Math. Phys. **219** (2001), no. 2, 271–322

[61] G. P. Scott and G. A. Swarup, *An algebraic annulus theorem,* Pacific J. Math. **196** (2000), no. 2, 461–506

[62] G. P. Scott and G. A. Swarup, *Canonical splittings of groups and 3-manifolds,* Trans. Amer. Math. Soc. **353** (2001), no. 12, 4973–5001

[63] Z. Sela, *Structure and rigidity in (Gromov) hyperbolic groups and discrete groups in rank* 1 *Lie groups. II,* Geom. Funct. Anal. **7** (1997), no. 3, 561–593

[64] Y. Shalom, *Random ergodic theorems, invariant means and unitary representation.* Lie groups and ergodic theory (Mumbai, 1996), 273–314, Tata Inst. Fund. Res. Stud. Math., **14**, Tata Inst. Fund. Res., Bombay, 1998

[65] Y. Shalom, *Expander graphs and amenable quotients.* Emerging applications of number theory (Minneapolis, MN, 1996), 571–581, IMA Vol. Math. Appl., **109**, Springer, New York, 1999

[66] H. Short, *Quasiconvexity and a theorem of Howson's,* in "Group theory from a geometrical viewpoint (Trieste, 1990)", 168–176, World Sci. Publishing, River Edge, NJ, 1991

[67] G. A. Swarup, *Geometric finiteness and rationality,* J. Pure Appl. Algebra **86** (1993), no. 3, 327–333

[68] G. A. Swarup, *Proof of a weak hyperbolization theorem,* Q. J. Math. **51** (2000), no. 4, 529–533

[69] W. Woess, *Random walks on infinite graphs and groups - a survey on selected topics,* Bull. London Math. Soc. **26** (1994), 1–60.

[70] W. Woess, *Random walks on infinite graphs and groups,* Cambridge Tracts in Mathematics, **138**. Cambridge University Press, Cambridge, 2000

Department of Mathematics, University of Illinois at Urbana-Champaign, 1409 West Green Street, Urbana, IL 61801, USA

*E-mail address*: kapovich@math.uiuc.edu